\documentstyle[amssymb]{article}
\newcommand{\PP}{{\mathbb {P}}}
\newcommand{\C}{{\mathbb {C}}}
\newcommand{\Q}{{$\mathbb {Q}$}}
\newcommand{\Z}{{\mathbb {Z}}}
\begin{document}

\title{K\"ahler-Einstein Metrics for some quasi smooth log del Pezzo Surfaces}
\author{Carolina Araujo}
\maketitle

\section{Introduction}

Recently Johnson and Koll\' ar determined the complete list of anticanonically
 embedded quasi smooth log del Pezzo surfaces in weighted projective 3-spaces 
[JK]. Then, using techniques developed earlier by Demailly and Koll\' ar [DK],
 they proved that many of those surfaces admit a K\" ahler-Einstein metric,
 and that some of them do not have tigers. Later Boyer, Galicki and Nakamaye
 showed the existence of K\" ahler-Einstein metric for $Z_{16} \subset 
\PP ^3$ (1,3,5,8) [BGN], one of the log del Pezzo surfaces for which 
the question was left open by Johnson and Koll\' ar.

The aim of this paper is to settle the question of the existence of
 K\"ahler-Einstein metrics and tigers for those surfaces in [JK] for which 
the question was still open. In order to show the existence of
 K\"ahler-Einstein metrics for these surfaces, we will use a criterion
 described in the work of Nadel [Na] and Demailly and Koll\'ar [DK].

\section{Anticanonically embedded quasi-smooth Log del Pezzo surfaces in weighted Projective 3-spaces}

For positive integers $q_i$ let $\PP (q_0,q_1,q_2,q_3)$ denote the
\emph{weighted projective 3-space} with weights $q_0 \leq q_1 \leq q_2 \leq
 q_3$. (See [Dol] or [Fle] for basic definitions and results). When there is
no ambiguity we abbreviate $\PP= \PP(q_0,q_1,q_2,q_3)$. We always assume
 that any 3 of the $q_i$ are relatively prime. Let $x_0, x_1, x_2, x_3$ denote
 the corresponding weighted projective coordinates. $\PP (q_0,q_1,q_2,q_3)$
 can be viewed as the quotient 
$$
\C^4\setminus\{0\}\  /\  \C^*(q_0,q_1,q_2,q_3)
$$
where $\C^*(q_0,q_1,q_2,q_3)$ denotes the $\C^*$-action on 
$\C^4\setminus\{0\}$ 
given by
$$
(z_0,z_1,z_2,z_3)\mapsto (t^{q_0}z_0,t^{q_1}z_1,t^{q_2}z_2,t^{q_3}z_3)
$$
The action of $\C^*$ on $\C^4 \setminus \{0\}$ induces the action 
$\Z_{q_i}(q_j,q_k,q_l)$ 
on $V_i=\{z_i=1\}\subset \C^4$ which, after identifying $V_i$ with $\C^3$, can 
be written as 
$$
(z_j,z_k,z_l)\mapsto (\xi ^{q_j}z_j,\xi ^{q_k}z_k,\xi ^{q_l}z_l)
$$
where $\xi$ is a primitive $q_i$th root of unity. Hence the affine chart 
$U_i\subset\PP$ \  where $x_i \neq 0$ can be identified with 
$$
\C^3/\Z_{q_i}(q_j,q_k,q_l)
$$
These are called the \emph{orbifold charts}. Write $p_i:\C^3 \to U_i$
for the natural projection.

$\PP (q_0,q_1,q_2,q_3)$ has an index $q_i$ singularity at 
$P_i=(x_j=0)_{j\neq i}$, and an index $(q_i,q_j)$ singularity along the line 
$(x_k=0)_{k\neq i,j}$.

For every $d\in \Z$ there is a rank 1 sheaf $\mathcal{O}$$_{\PP}(d)$, and the 
sections of $\mathcal{O}$$_{\PP}(d)$ consist of the homogeneous polynomials 
$f(x_0,x_1,x_2,x_3)$ of weighted degree d.

Let $X=X_d$ be a member of $|\mathcal{O}$$_{\PP}(d)|$. If $X$ does not contain 
any of the singular lines of $\PP$ , then the adjunction formula 
$$
K_X \cong {\mathcal{O}} _{\PP}(K_{\PP}+X)|_X \cong 
{\mathcal{O}} _{\PP}(d-(q_0+q_1+q_2+q_3))|_X
$$
holds.

A \emph{log del Pezzo} surface is a projective surface with quotient 
singularities such that its anticanonical class is ample.

If $X$ does not contain any of the singular lines, then, by the adjunction
formula, its anticanonical class is ample if and only if $d<q_0+q_1+q_2+q_3$. 
Here we are interested in the case when $d$ is as large as possible, i.e. 
$d=q_0+q_1+q_2+q_3-1$. In this case we say that $X$ is 
\emph{anticanonically embedded} in $\PP (q_0,q_1,q_2,q_3)$.

We say that $X$ is \emph{quasi smooth} if the pre image of $X$ under the 
quotient map
$$
\C^4 \setminus \{0\} \to \C^4 \setminus \{0\}\ /\ \C^*(q_0,q_1,q_2,q_3)=
\PP(q_0,q_1,q_2,q_3)
$$ 
is smooth.

We write conditions on the weights $q_0,q_1,q_2,q_3$ for the general surface 
$X_d \subset \PP (q_0,q_1,q_2,q_3)$ of degree $d=q_0+q_1+q_2+q_3-1$ to be 
quasi smooth, and for the adjunction formula to hold. These conditions then
imply that the general $X_d$ is an \emph{anticanonically embedded quasi smooth
log del Pezzo} surface: 

\paragraph{Conditions 2.1}{\it{For $X_d$ to be quasi smooth (see [Fle], 8.5):

(I) For every $i$ there is a $j$ and a monomial $x_i^{m_i}x_j$ of degree $d$.

(II) For all distinct $i$ and $j$, either there is a monomial 
$x_i^{b_i}x_j^{b_j}$ of degree $d$, or there are monomials 
$x_i^{c_i}x_j^{c_j}x_k$ and $x_i^{d_i}x_j^{d_j}x_l$ of degree $d$, 
with $k$ and $l$ distinct.

For $X_d$ not to contain any of the singular lines of $\PP$ 
\  (which implies that the adjunction formula holds):

(III)  For all distinct $i$ and $j$ such that $(q_i,q_j)>1$, there is a 
monomial $x_i^{b_i}x_j^{b_j}$ of degree $d$.}}  $\Box$

\section{Tigers and the Existence of K\"ahler-Einstein Metrics}

We start this section by giving some definitions from the log category.
We refer to [KM] for a detailed introduction.

\paragraph{Definition}Let $X$ be a normal surface, and $D$ a \Q-divisor on 
$X$ such that $m(K_X+D)$ is Cartier for some $m>0$. Let $f:Y \to X$ be a proper
birational morphism from a smooth surface $Y$. Then there is  unique \Q-divisor
$\sum e_iE_i$ on $Y$ such that
$$
K_Y \equiv f^*(K_X+D)+\sum e_iE_i \ \ and \ \ \  f_*\sum e_iE_i =-D.
$$
We say that the pair $(X,D)$ is \emph{canonical} (resp. \emph{klt}, resp.
\emph{log canonical}) if $e_i \ge 0$ (resp. $e_i>-1$, resp. $e_i \ge -1$) for 
every $f$ and every $i$. 

\paragraph{Remark}If $(X,D)$ and $(X,D')$ are both canonical (resp.
klt, resp. log canonical), then so is $(X,\alpha D+(1-\alpha )D')$ 
for any $0\leq \alpha \leq 1$. This is a very easy result, but it will
be very useful in our applications.

\paragraph{Definition}[KMcK] Let $X$ be a normal surface. A 
\emph{tiger} on $X$ is an effective \Q-divisor $D\equiv -K_X$ 
such that $(X,D)$ is not klt. \bigskip

We use the following sufficient condition for the existence of 
K\"ahler-Einstein metric on log del Pezzo surfaces:

\paragraph{Theorem 3.1}([Na],[DK]) {\it{Let $X$ be a log del Pezzo surface. 
If there is an $\epsilon >0$ such that $(X,\frac{2+\epsilon}{3}D)$ is klt 
for every effective \Q-divisor $D \equiv -K_X$, then X admits a 
K\"ahler-Einstein metric.}} $\Box$ \bigskip

In order to prove the existence of K\"ahler-Einsten metric (resp. the 
nonexistence of tigers) for a surface $X_d \subset\PP(q_0,q_1,q_2,q_3)$, 
we will show that, for all effective \Q-divisors $D \equiv 
-\frac{2+ \epsilon}{3}K_X$ (resp $D \equiv K_X$), $(X,D)$ is not klt. 
The definition of a klt pair 
$(X,D)$, however, involves understanding all resolutions of singularities 
of $(X,D)$. Instead, we will use the following multiplicity conditions:

\paragraph{Conditions 3.2}[JK] {\it{Let $X=X_d$ be a hypersurface of 
weighted degree $d$ in $\PP(q_0,q_1,q_2,q_3)$, and assume that $X$ is 
quasi-smooth and has only isolated (quotient) singularities. Write
$$
p_i\ :\ \C^3\ \to \ U_i\ =\ (x_i\neq0)\ \cong\ 
\C^3\ /\ \Z_{q_i}(q_j,q_k,q_l)
$$
for the natural projections to the orbifold charts.
For an effective \Q-divisor $D$ on 
$X$, $(X,D)$ is klt provided that the following 3 conditions hold:

(0) $D$ does not contain any irreducible component with coefficient $\ge 1$.

(1) For every smooth point $P \in X$, $mult_P D \leq 1$

(2) For every singular point $P \in X \cap U_i$, $mult_Q p_i^*D \leq 1$, 
where $p_i(Q)=P$.}} $\Box$ \bigskip

In our applications, condition (3.2.0) can be verified right away. To check 
the other 2 conditions, we will bound the multiplicities above by suitable 
intersection numbers.

\subsection{The Smooth Points}

Let $P \in X^0$, the nonsingular locus of X. Suppose we can find a positive 
integer $l$ such that, for any effective \Q-divisor
$D \equiv -K_X$, there is a divisor $F \in 
|{\mathcal{O}}_X(l)|$ such that:

\noindent (1) $F=F'+B$, where $F'$ is an effective Weil divisor, and $B$
is nef.

\noindent (2) $F'$ passes through $P$.

\noindent (3) $F'$ does not contain any irreducible component of $D$.

\noindent Then we can intersect $D$ with $F'$, and obtain:
$$
mult_P D \leq (D\cdot F') \leq (D\cdot F)
= \frac{ld}{q_0q_1q_2q_3}
$$ 
We would like this bound to be $\leq 1$ in order to verify condition (3.2.1).
Unfortunately in most of the cases we can not choose an $l$ that works for 
all $P \in X^0$. Instead, we will find an $l$ that works for a dense 
open set of $X^0$, say $X^0 \setminus C$, and then show by other methods
that $mult_P D \leq 1$ for $P \in C \cap X^0$.

\paragraph{Lemma 3.3}{\it{Let $\PP = \PP (q_0,\ldots ,q_n)$ denote the 
weighted projectice n-space with weights $q_0,\ldots ,q_n$. 
Let $l$ be a positive integer such that for every $i \neq 0$ there
are at least 2 distinct monomials in $H^0(\PP, {\mathcal{O}}_{\PP}(l))$
of the form $x_0^{\alpha}x_i^{\beta}$. 
Then, for every $P \in \PP \setminus (x_0=0)$, and every divisor $D$ on 
$\PP$, there is a divisor $F\in |{\mathcal{O}}_{\PP}(l)|$ that can be 
written as $F=F'+a(x_0=0)$, where $a$ is a non negative integer, and $F'$
is an effective Weil divisor passing through 
$P$ and not containing any irreducible component of $D$.}} \bigskip

\noindent Proof. We choose from each irreducible component $D_j$ of $D$,
$D_j$ not in $(x_0=0)$,
a point $Q_j \in D_j \setminus (x_0=0)$
(later we will impose a further condition on the $Q_j$).
We shall prove that there is a polynomial $F_j \in 
H^0(\PP, {\mathcal{O}}_{\PP}(l))$
that vanishes at $P$ but not at $Q_j$. We then
choose $F$ to be a suitable linear combination of the $F_j$.

Let $Q=Q_j$. Fix representations 
$P=(1:a_1:\ldots:a_n)$ and $Q=(1:b_1:\dots:b_n)$ in weighted projective
coordinates (Notice that these representations are not unique). If 
$a_i=0\neq b_i$ or $a_i\neq 0=b_i$ for some $i$, then it is easy to
find a polynomial in $H^0(\PP, {\mathcal{O}}_{\PP}(l))$
that vanishes at $P$ but not at $Q$. So we assume that for every $i$, 
either $a_i=0=b_i$ or $a_i\neq 0\neq b_i$.

For each $i \ne 0$ fix $x_0^{\alpha_i}x_i^{\beta_i}$ and 
$x_0^{\gamma_i}x_i^{\delta_i}$ distinct monomials in 
$H^0(\PP, {\mathcal{O}}_{\PP}(l))$, and assume $\delta_i>\beta_i$. Define:
$$
H_i\ =\ a_i^{\delta_i}x_0^{\alpha_i}x_i^{\beta_i}\ - \ 
a_i^{\beta_i}x_0^{\gamma_i}x_i^{\delta_i}\ 
\in H^0(\PP, {\mathcal{O}}_{\PP}(l))
$$
Then $H_i(P)=0$, and $H_i(Q)=0$ if and only if 
$b_i=\zeta_ia_i$, where $\zeta_i^{\delta_i-\beta_i}=1$. 

Therefore, if $Q$ is outside the finite set $\{(1:\zeta_1a_1:\ldots :
\zeta_na_n)\ |\ \zeta_i^{\delta_i-\beta_i}=1\}$
(and we can certainly impose that condition when we choose the $Q_j$
above), we get a polynomial in $H^0(\PP, {\mathcal{O}}_{\PP}(l))$
that vanishes at $P$ but not at $Q$. 
$\Box$
\bigskip

Now we will use Lemma 3.3 to find a linear system  
$|{\mathcal{O}}_{X_d}(l)|$ on $X_d$ that will give us a bound for 
$mult_PD$, at least for $P\in X_d\cap (x_0\neq 0)$.

\paragraph{Corollary 3.4}{\it{Let $X=X_d \subset \PP(q_0,q_1,q_2,q_3)$ 
be an anticanonically embedded quasi smooth log del Pezzo surface. 
Let $\pi_3:X \to \PP=\PP(q_0,q_1,q_2)$ denote the projection from 
$P_3=(0,0,0,1)$. Assume $\pi_3$ has only finite fibers. Let $l$ be a 
positive integer such that, for $i=1,2$, there are at least 2 monomials 
in $H^0(\PP,{\mathcal{O}}_{\PP}(l))$ of the form 
$x_0^{\alpha}x_i^{\beta}$.

Then, for every $P\in X^0\setminus (x_0=0)$ and every effective \Q-divisor 
$D\equiv -K_X$, there is a divisor $F\in |{\mathcal{O}}_X(l)|$ 
that can be written as $F=F'+a(x_0=0)$, where $a$ is a non negative 
integer, and $F'$ is an effective Weil divisor passing 
through  $P$ and not containing any irreducible component of $D$. Hence
$mult_P D\leq \frac{ld}{q_0q_1q_2q_3}$.}}\bigskip

\noindent Proof. Set $P'=\pi_3(P)\in \PP$, and $D'=\pi_{3*}D$.
 By the lemma, there is a divisor $E \in 
|{\mathcal{O}}_{\PP}(l)|$ that can be writen as $E=E'+a(x_0=0)$, 
where $a$ is a non negative integer, and $E'$ is an effective Weil divisor
passing through $P'$ and not containing any
irreducible component of $D'$. Since $\pi_3$ has only finite fibers, 
we can take $F$ to be
$\pi_3^*E\in \pi_3^*|{\mathcal{O}}_{\PP}(l)| \subset |{\mathcal{O}}_X(l)|$.
$\Box$
\bigskip

Now we have to deal with the smooth points of $X$ in $(x_0=0)$.

\paragraph{Lemma 3.5}{\it{Let $X=X_d \subset \PP(q_0,q_1,q_2,q_3)$ 
be an anticanonically embedded quasi smooth log del Pezzo surface. 
Assume that $C=X\cap (x_0=0)$ is irreducible and smooth outside
the singular locus of X. If $d\leq q_1q_2q_3$, then for every effective
\Q-divisor
$D\equiv -K_X$ and every $P\in X^0\cap C$ we have $mult_P D\leq 1$.}}
\bigskip

\noindent Proof. Fix $P\in X^0\cap C$. Since $C$ is smooth outside 
the singular locus of $X$, $mult_P C=1$. Given $D\equiv -K_X$, write
$D=\alpha C+(1-\alpha q_0)D'$, where $0\leq \alpha \leq 1/q_0$,
$D' \equiv -K_X$, and $C$ is not contained in the support of $D'$. Then
$$
mult_P D'\ \leq\ (C\cdot D')\ \leq \ \frac{d}{q_1q_2q_3}\ \leq\ 1
$$
And hence $mult_P D\leq \alpha +(1-q_0\alpha)\leq 1$.
$\Box$
\bigskip
 
Now we put these results together to obtain an arithmetical 
condition on the weights $q_i$ for every effective \Q-divisor 
$D\equiv -K_{X_d}$ 
to have multiplicity $\leq 1$ at the smooth points of $X_d$. 

\paragraph{Lemma 3.6}{\it{Let $X=X_d \subset \PP(q_0,q_1,q_2,q_3)$  
be an anticanonically embedded quasi smooth log del Pezzo surface. Assume 
that:

\noindent (1) The curve $C=X\cap (x_0=0)$ is irreducible and smooth
outside the singular locus of $X$.

\noindent (2) There is a monomial of the form $x_3^{\alpha}$ appearing
in the equation of $X$ with nonzero coefficient. 

\noindent (3) $d\leq q_1q_3$.

Then $mult_P D\leq 1$ for every smooth point $P\in X$ and every 
\Q-divisor $D\equiv -K_X$.}}\bigskip

\noindent Proof. (1), (3) and Lemma 3.5 together give the result when 
$P\in C$. (2) implies that the projection $\pi_3$
has only finite fibers. We then apply Lemma 3.4 with $l=q_0q_2$, 
and obtain the result for $P\not\in C$.
$\Box$

\paragraph{Remarks}

\noindent (i) In each of our applications, conditions (3.6.1) and 
(3.6.2) above will be forced by conditions (2.1).

\noindent (ii) The same result holds if we replace $x_3^{\alpha}$ by
$x_1^{\alpha}$ or $x_2^{\alpha}$ in (3.6.1), and the condition 
$d\leq q_1q_3$ by $d\leq q_1q_2$ in (3.6.2).

\subsection{The Singular Points} 

\paragraph{Estimate 3.7}{\it{Let $P\in X \cap U_i$ be a singular point, 
and let $p_i:(S,Q)\to (X,P)$ be a local orbifold chart. If 
$(x_j=x_k=0)\not\subset X$ then, for every effective \Q-divisor 
$D\equiv -K_X$,
$$
mult_Qp_i^*D \ \ \leq \ \ \frac{d}{min(q_j,q_k)\cdot q_l}
$$
}}

\noindent Proof. Consider the linear system
$$
|p_i^*(x_j^{q_k}),p_i^*(x_k^{q_j})| \ \ \subset \ \  
p_i^*|{\mathcal{O}}_X(q_jq_k)|
$$
By hypothesis this linear system does not have any fixed component. So we 
can intersect a general member of it with $p_i^* D$, obtaining:
$$
mult_Qp_i^*D \ \leq \frac{q_i}{min(q_j,q_k)}(D \cdot 
{\mathcal{O}}_X(q_jq_k))\ = \frac{d}{min(q_j,q_k)\cdot q_l}
$$
$\Box$
\medskip

When this bound is $>1$ but $\leq 2$, we can look closer at what happens 
near $Q$, that is, we can blow up this point and then use Shokurov's 
Inversion of Adjunction:

\paragraph{Theorem 3.8}([KM], 5.50): {\it{Let $S$ be a smooth surface, 
$Q \in S$, and $B$ be an effectice \Q-divisor such that $mult_QB \leq 2$. Let 
$\pi : S' \to S$ be the blow up of Q with exceptional divisor $E$. If 
$(\pi ^{-1}_* (B))|_E$ is a sum of points, all with coefficients $\leq 1$,
then $(S,B)$ is log canonical near $Q$.}} $\Box$

\section{The 6 Surfaces}
 
\paragraph{Theorem 4.1}{\it{Let $X_d \subset \PP(q_0,q_1,q_2,q_3)$ 
denote any anticanonically embedded quasi smooth log del Pezzo surface 
in the weighted projective 3-space with weights $q_0,q_1,q_2,q_3$. 
Then:

\noindent (1) $X_{36}\subset \PP(3,5,11,18)$ does not have a tiger, and 
hence admits a K\" ahler-Einstein metric.

\noindent (2) $X_{28}\subset \PP(3,5,7,14)$ does not have a tiger, and 
hence admits a K\" ahler-Einstein metric.

\noindent (3) $X_{25}\subset \PP(3,5,7,11)$ does not have a tiger, and 
hence admits a K\" ahler-Einstein metric.

\noindent (4) $X_{18}\subset \PP(2,3,5,9)$ does not have a tiger, and 
hence admits a K\" ahler-Einstein metric.

\noindent (5) $X_{15}\subset \PP(1,3,5,7)$ admits a K\" ahler-Einstein 
metric provided that the coefficient of the monomial $x_1x_2x_3$ in the 
equation of $X_{15}$ is nonzero.

\noindent (6) $X_{10}\subset \PP(1,2,3,5)$ admits a K\" ahler-Einstein 
metric.}}

\paragraph{Remarks}

\noindent (i) The existence of K\"ahler-Einstein metric for the surfaces
in (4.1.1), (4.1.2) and (4.1.3) was already known.

\noindent (ii) $(x_0=0)$ is a tiger on the surfaces in 
(4.1.5) and (4.1.6).

\noindent (iii) In (4.1.5), if the coefficient of $x_1x_2x_3$ in the 
equation of $X_{15}$ is zero, then the criterion of Theorem 3.1 does not 
apply. Indeed, $c \cdot  (x_0=0)$ is not klt for any $c>8/15$.
However, it is possible that even in this case $X_{15}$ admits a 
K\"ahler-Einstein metric.\bigskip

\noindent Proof. (1) Lemma 3.6 and Estimate 3.7 imply that  
Conditions 3.2 are verified.\medskip

\noindent (2) Corollary 3.4, 
with $l=21$, implies that for every smooth point $P\in X\setminus 
(x_0=0)$, and every effective \Q-divisor $D\equiv -K_X$,
$mult_P D\leq 2/5$. The curve $X\cap (x_0=0)$ has 2 irreducible 
components, so we cannot apply Lemma 3.6. Instead we apply Corollary 3.4
again with $x_0$ replaced by $x_2$, and $l=35$. This shows that for
every smooth point $P\in X\setminus (x_2=0)$, and every effective 
\Q-divisor $D\equiv -K_X$, $mult_P D\leq 2/3$. Estimate 3.7 takes care
of the singular points. \medskip

\noindent (3) Lemma 3.6 takes care of the smooth points of $X$. 
The singular points
of $X$ are $P_0=(1:0:0:0)$, $P_2=(0:0:1:0)$ and $P_3=(0:0:0:1)$. Estimate
3.7 takes care of $P_0$ and $P_2$, but only gives that 
$mult_{Q_3}p_3^*D\leq 25/21$ for any effective \Q-divisor 
$D\equiv -K_X$. But we can improve this bound:

We first notice that $C=X\cap (x_0=0)\in |{\mathcal{O}}_X(3)|$ is 
irreducible, $p_3^*C \cong (ay_1^5+by_2^2=0)\subset \C^2(y_1,y_2)$,
and so $mult_{Q_3}p_3^*C=2$ (Here $a$ and $b$ are the coefficients 
of $x_1^5$ and $x_2^2x_3$ in the equation of $X$, which are nonzero 
by condition (2.1.I)).
 Write $D=\alpha C+(1-3\alpha)D'$, where
$0\leq \alpha \leq 1/3$, $D'\equiv -K_X$, and $C$ is not contained 
in the support of $D'$. Then:
$$
mult_{Q_3}p_3^*D'\leq \frac{1}{2} (p_3^*C \cdot p_3^*D')\leq
\frac{11}{2} ({\mathcal{O}}_X(3)\cdot {\mathcal{O}}_X(1))=\frac{5}{14}
$$
So $mult_{Q_3} p_3^*D\leq 2/3$, and we have taken care of all singular 
points.\medskip

\noindent (4) Lemma 3.6 and Estimate 3.7 imply that Conditions 3.2 are 
verified.\medskip

\noindent (5) Let $a$ and $b$ be the coefficients of $x_1^5$ 
and $x_2^3$ in the
equation of $X$. Condition (2.1.I) implies that $a$ and $b$ are nonzero.
If the coefficient of $x_1x_2x_3$ in the equation of $X$ is zero,
then $C=X\cap (x_0=0)\cong (ax_1^5+bx_2^3=0) \subset \PP(3,5,7)$ is
$\equiv -K_X$, and $(X,\frac{2}{3}C)$ is not klt. (To see this, notice that 
$(S=p_3^{-1}(X),p_3^*(\frac{2}{3}C))$ 
is not klt by ([Kol], 8.15), and then apply 
([KM], 5.20)). In this case the criterion of Theorem 3.1 can not be applied.
From now on we assume that the coefficient of $x_1x_2x_3$ in the equation 
of $X$ is nonzero. 

Let $D\equiv -K_X$ be an effective \Q-divisor. We shall show that 
$(X,D)$ is log canonical. This is more than enough to apply Theorem
3.1. 

Write $D=\alpha C+(1-\alpha) D'$, with $0\leq \alpha \leq 1$, 
$D\equiv -K_X$, and $C$ not contained in the support of $D'$. It
is enough to show that both $(X,C)$ and $(X,D')$ are log canonical.

Lemma 3.6 shows that any effective \Q-divisor $D\equiv -K_X$
has multiplicity $\leq 1$ at the smooth
points of $X$, and hence $(X,D)$ is log canonical at the smooth points 
of $X$. 

$P_3=(0:0:0:1)$ is the only singular point of $X$.
Notice that $p_3^*C$ has a simple node at $Q_3$. In particular 
$mult_{Q_3}p_3^*C= 2$, and we can use Shokurov's inversion of 
adjunction. Let the notation be as in Theorem 3.8. Then 
$(\pi^{-1}_*(p_3^*C))|_E$ is a sum of 2 points, each with coefficient $1$.
We conclude that $(X,C)$ is log canonical near $P_3$.

Now for $D'$ we have:
$$
mult_{Q_3}p_3^*D'\ \leq\ \frac{1}{2}(p_3^*D'\cdot p_3^*C)\ \leq\ 
\frac{7}{2}(D'\cdot C)\ =\frac{1}{2}
$$
and hence $(X,D')$ is log canonical near $P_3$. 
\medskip

\noindent (6) Let $C=X\cap (x_0=0)\cong 
(ax_1^5+bx_1^2x_2^2+cx_1x_2x_3+dx_3^2=0)=
(ax_1^5+(\alpha x_1x_2+\beta x_3)(\gamma x_1x_2+\delta x_3)=0)\subset 
\PP(2,3,5)$. Here $a$, $b$, $c$, and $d$ are the coefficients of 
$x_1^5$, $x_1^2x_2^2$, $x_1x_2x_3$ and $x_3^2$ in the equation of 
$X$, and $a,d\neq 0$ by condition (2.1.I).

Let $D\equiv -K_X$ be an effective \Q-divisor. We shall show that 
$(X,\frac{7}{10}D)$ is log canonical. This is more than enough to 
apply Theorem 3.1. 

Write $D=\lambda C+(1-\lambda) D'$, with $0\leq \lambda \leq 1$, 
$D'\equiv -K_X$, and $C$ not contained in the support of $D'$. It
is enough to show that both $(X,\frac{7}{10}C)$ and $(X,\frac{7}{10}D')$ 
are log canonical.

Lemma 3.6 shows that any effective \Q-divisor $D\equiv -K_X$
has multiplicity $\leq 1$ at the smooth
points of $X$, and hence $(X,\frac{7}{10}D)$ is log canonical at 
the smooth points of $X$. 

$P_2=(0:0:1:0)$ is the only singular point of $X$.
$mult_{Q_2}p_2^*C= 2$, and we can use Shokurov's inversion of 
adjunction. Let the notation be as in Theorem 3.8.
We have 2 cases: If $(\alpha:\beta)\neq
(\gamma:\delta)$, then $p_2^*C$ has a simple node at $Q_2$, and 
$(\pi^{-1}_*(p_2^*C))|E$ is a sum of 2 points, each with coefficient 1.
In this case $(X,C)$ is log canonical near $P_2$. If $(\alpha:\beta)=
(\gamma:\delta)$, then $p_2^*C$ has only one tangent direction at $Q_2$,
and $(\pi^{-1}_*(p_2^*C))|E$ is a single point, with coefficient 2. In this 
case $(X,C)$ is not log canonical near $Q_2$. However, using ([Kol],8.15)
and ([KM],5.20), we conclude that $(X,\frac{7}{10}C)$ is log canonical 
near $Q_2$.

Now for $D'$ we have:
$$
mult_{Q_2}p_2^*D'\ \leq\ \frac{1}{2}(p_2^*D'\cdot p_2^*C)\ \leq\ 
\frac{3}{2}(D'\cdot C)\ =\frac{1}{2}
$$
and hence $(X,\frac{7}{10}D')$ is log canonical near $P_2$. 
\ $\Box$

\paragraph{Acknowledgements} I would like to thank J. Koll\'ar for very 
useful comments and corrections. I would also like to thank the R\'enyi 
Institute of Mathematics for the hospitality during the conference on
Higher Dimensional Varieties and Rational Points - September 2001, when 
part of this paper was written. Partial financial support was provided
by the CNPq (Conselho Nacional de Desenvolvimento Cient\'ifico e 
Tecnol\'ogico - Brazil).

\bigskip

\end{document}